\documentclass{amsart}

\usepackage{amsmath,amssymb,amsthm,latexsym}

\newtheorem{theorem}{Theorem}
\newcommand{\bt}{\begin{theorem}}
\newcommand{\et}{\end{theorem}} 

\newtheorem{lemma}{Lemma}
\newcommand{\bl}{\begin{lemma}}
\newcommand{\el}{\end{lemma}}

\newtheorem{corollary}{Corollary}
\newcommand{\bc}{\begin{corollary}}
\newcommand{\ec}{\end{corollary}} 

\newtheorem{problem}{Problem}
\newcommand{\bprob}{\begin{problem}}
\newcommand{\eprob}{\end{problem}}

\newcommand{\beq}{\begin{equation}}
\newcommand{\eeq}{\end{equation}}
\newcommand{\N}{\ensuremath{ \mathbf N }}
\newcommand{\Z}{\ensuremath{\mathbf Z}}

\newcommand{\R}{\ensuremath{\mathbf R}}
\newcommand{\benum}{\begin{enumerate}}
\newcommand{\eenum}{\end{enumerate}}

\newcommand{\mcf}{\mathcal{F}}
\newcommand{\mch}{\mathcal{H}}
\newcommand{\mbh}{\mathbf{H}}

\DeclareMathOperator{\qqand}{\qquad\text{and}\qquad}
\DeclareMathOperator{\card}{\text{card}}

\title[Intersections of  product sets and sumsets]
{Problems and results on intersections of product sets and sumsets in semigroups}

\author{Melvyn B.  Nathanson}
\address{Department of Mathematics\\Lehman College (CUNY)\\Bronx, NY 10468}
\email{melvyn.nathanson@lehman.cuny.edu}

\date{\today}

\thanks{Supported in part by  PSC-CUNY Research Award Program grant 66197-00 54.}

\begin{document}

\maketitle
\begin{abstract}
For every subset $A$ of a semigroup $S$, let $A^h$  be the set of all products 
of $h$ elements of $S$. If $(A)_{q\in Q}$ is a family of subsets of $S$, then 
$A = \bigcap_{q \in Q} A_q$ satisfies $A^h \subseteq  \bigcap_{q \in Q} A_q^h$.  
The product intersection set 
$H(A_q) = \left\{h \in \mathbf{N}: A^h = \bigcap_{q \in Q} A_q^h \right\}$  
is investigated. 
\end{abstract}

\section{The product  and sum  intersection sets}

Let $\N = \{1,2,3,\ldots\}$ be the set of positive integers.  
The letter $h$ will always denote a positive integer.

Let $S$ be a semigroup, written multiplicatively and not necessarily abelian, 
with subsets $A_1,\ldots, A_h$.  
The \emph{product set} of these subsets is the set 
\[
A_1A_2\cdots A_h = \{a_1a_2\cdots a_h: a_i \in A_i \text{ for all } i = 1,\ldots, h\}.
\]
If $A_i = A$ for all $i=1,\ldots, h$, then the \emph{$h$-fold product set} is 
\[
A^h = \{a_1a_2\cdots a_h: a_i \in A \text{ for all } i = 1,\ldots, h\}.
\]
If the semigroup $S$ is abelian and written additively, then we have the \emph{sumset} 
\[
A_1 + A_2 + \cdots + A_h = \{a_1 + a_2 + \cdots + a_h: a_i \in A_i \text{ for all } i = 1,\ldots, h\} 
\]
and  the \emph{$h$-fold sumset}
\[
hA = \{a_1 + a_2 + \cdots + a_h: a_i \in A \text{ for all } i = 1,\ldots, h\}.
\]

Let $Q$ be an infinite set and let $(A_q)_{q\in Q}$ be a family 
of subsets of the semigroup $S$ such the the set of sets $\{A_q:q\in Q\}$ is infinite.  
The set $Q$ is called an \emph{index set}. 
Let $\mcf_Q(S)$ be the set of all infinite families $(A_q)_{q\in Q}$ 
of subsets of $S$ indexed by $Q$.  
For all $(A_q)_{q\in Q} \in \mcf_Q(S)$, we define 
\[
A = \bigcap_{q\in Q} A_q.
\]
For multiplicative semigroups, we have the multiplicative inclusions 
$A^h \subseteq A_q^h$   for all $q\in Q$ and $h \in \N$ and so  
\beq               \label{intersect:multiplicative}
A^h \subseteq \bigcap_{q\in Q} A_q^h.
\eeq
For additive abelian semigroups,  we have  the additive inclusions $hA\subseteq hA_q$ 
for all $q\in Q$ and $h \in \N$ and so  
\beq               \label{intersect:additive}
hA \subseteq \bigcap_{q\in Q} hA_q. 
\eeq
For some families $(A_q)_{q \in Q}$, the 
relations~\eqref{intersect:multiplicative} and~\eqref{intersect:additive}  
are equalities, and for other sequences the inclusions are proper inclusions. 
For a multiplicative semigroup $S$ and index set $Q$, 
we define the \emph{product intersection set} 
\[
H(A_q) = H_Q^S(A_q) = \left\{ h \in \N: A^h  = \bigcap_{q \in Q} A_q^h \right\}. 
\] 
For an additive semigroup $S$ and index set $Q$,  
we define the \emph{sum intersection set} 
\[
 H(A_q) = H_Q^S(A_q) = \left\{ h \in \N: hA  = \bigcap_{q \in Q} hA_q \right\}. 
\] 
We have $1 \in H(A_q)$ for all $(A_q)_{q\in Q} \in \mcf_Q(S)$.  

There are many unsolved problems about the intersections 
of $h$-fold product sets of subsets 
of a semigroup $S$ with respect to an index set $Q$.

\bprob              \label{interset:problem:H}
For a given family $(A_q)_{q\in Q} \in \mcf_Q(S)$, compute $H(A_q)$. 
\eprob 

\bprob                 \label{interset:problem:mathcal-H}
For a given subset $A$ of $S$, determine the set 
\[
\mch_Q(A) = \left\{ H(A_q) : (A_q)_{q\in Q} \in \mcf_Q(S) \text{ and } A = \bigcap_{q\in Q} A_q \right\}.
\] 
\eprob 

For example, the problem $\mch_Q(\emptyset)$ is to determine all subsets $X$ of \N\ such that 
there exists a family $(A_q)_{q\in Q} \in \mcf_Q(S)$ with $\bigcap_{q\in Q} A_q = \emptyset$ 
such that $\bigcap_{q\in Q} A_q^h = \emptyset$ if and only if $h \in X$.

\bprob                  \label{interset:problem:mathbf-H}
Determine the set 
\[
\mathbf{H}_Q(S)  = \bigcup_{A\subseteq S}  \mch_Q(A) 
 = \left\{  H(A_q):    (A_q)_{q\in Q} \in \mcf_Q(S) \right\}. 
\]
\eprob

The set $\mathbf{H}_Q(S)$ is the ``master set'' of all product intersection sets $H(A_q)$ 
for all families $(A_q)_{q\in Q} \in \mcf_Q(S)$.   
The following is a special case of Problem~\ref{interset:problem:mathbf-H}.  

\bprob                                 \label{interset:problem:noX}
Is there a set $X$ of positive integers with $1 \in X$ such that $X \notin \mathbf{H}_Q(S_0)$ 
for the semigroup $S_0$? 

Is there a set $X$ of positive integers with $1 \in X$ such that $X \notin \mathbf{H}_Q(S)$ 
for all semigroups $S$? 

\eprob

Problem~\ref{interset:problem:noX} asks if there exists a set $X$ of positive integers with $1 \in X$ 
such that there exists no family $(A_q)_{q=1}^{\infty} \in \mcf_Q(S)$ 
with $\bigcap_{q\in Q}A_q = A$ such that $A^h = \bigcap_{q\in Q} A_q^h$ if and only if $h \in X$. 
 
Replacing $h$-fold product sets with $h$-fold sumsets, we have the analogous problems 
for additive abelian semigroups. 
Product intersection sets and sum intersection sets are of special interest in number theory 
when investigated for the additive and multiplicative semigroups \N\ and \Z. 

The purpose of this paper is to prove basic results about intersection sets and to construct 
explicit examples of families of subsets of a semigroup with interesting intersection properties. 

\section{Intersections of sequences of sets}

We often consider intersection sets with respect to  the index set $Q = \N$ 
and sequences of subsets $(A_q)_{q=1}^{\infty}$.  
A  sequence  $(A_q)_{q=1}^{\infty}$  of subsets of a set $S$ is \emph{decreasing}
 if $A_q \supseteq A_{q+1}$  for all $q \in \N$ and  
\emph{strictly decreasing}\index{decreasing sequence!strictly} 
if $A_q \supseteq A_{q+1}$  and $ A_q \neq A_{q+1}$ for all $q \in \N$.  
  A decreasing sequence   $(A_q)_{q=1}^{\infty}$  is 
\emph{asymptotically strictly decreasing}\index{sequence!asymptotically strictly decreasing} 
 if  $ A_q \neq A_{q+1}$ for infinitely many $q \in \N$. 
A decreasing sequence is asymptotically strictly decreasing if and only if it is not eventually constant.

 Let $\mcf_{\N}^*(S)$  be the set of all asymptotically strictly decreasing sequences 
$(A_q)_{q=1}^{\infty}$ of subsets of $S$.

\bprob                 \label{interset:problem:mathcal-HN}
For a given subset $A$ of $S$, determine the set 
\[
\mch_{\N}^*(A) = \left\{ H(A_q) : (A_q)_{q=1}^{\infty} \in \mcf_{\N}^*(S) \text{ and } A = \bigcap_{q=1}^{\infty} A_q \right\}.
\] 
\eprob

\bprob                  \label{interset:problem:mathbf-HN}
Determine the set 
\[
\mathbf{H}^*_{\N}(S)  = \bigcup_{A\subseteq S} \mch_{\N}^*(A) 
 = \left\{  H(A_q):   (A_q)_{q=1}^{\infty} \in \mcf_{\N}^*(S) \right\}. 
\]
\eprob

\bprob                                 \label{interset:problem:noXN} 
Is there a set $X$ of positive integers with $1 \in X$ such that $X \notin \mathbf{H}^*_{\N}(S_0)$ 
for the semigroup $S_0$? 

Is there a set $X$ of positive integers with $1 \in X$ such that $X \notin \mathbf{H}^*_{\N}(S)$ 
for all semigroups $S$? 
\eprob

 A subset $A$ of a set $S$ is  \emph{co-finite} if the complement of $A$ in $S$ is finite 
 and \emph{co-infinite} if the complement of $A$ in $S$ is infinite.  
 
 If  the sequence $(A_q)_{q=1}^{\infty}$  is strictly decreasing and $A =  \bigcap_{q=1}^{\infty} A_q$, 
 then, for all $q \in \N$, we have 
 \[
 S\setminus A \supseteq A_1\setminus A \supseteq A_1\setminus A_{q+1} 
 =  \bigcup_{r=1}^q A_r\setminus A_{r+1}.
 \]
Because the sets $A_r\setminus A_{r+1}$ are pairwise disjoint and nonempty, we have 
 \[
 \left| S \setminus A \right| \geq  \left|  \bigcup_{r=1}^q A_r\setminus A_{r+1}  \right| 
 = \sum_{r=1}^q   \left| A_r\setminus A_{r+1}  \right| \geq q 
 \]
for all $q \in \N$ and so $A$ is a co-infinite subset of $S$. 
Every asymptotically strictly decreasing sequence $(A_q)_{q=1}^{\infty}$ contains a strictly decreasing subsequence and so $A =  \bigcap_{q=1}^{\infty} A_q$ is also a co-infinite subset of $S$.

For many families  $(A_q)_{q \in Q} \in \mcf_{Q}(S)$, we have $H(A_q) = \{1\}$ 
or $H(A_q) = \N$.  
For example, if $S = \N$ and $Q= \N$, then the sequence of sets 
$(A^{\sharp}_q)_{q=1}^{\infty}$ defined by 
\[
A^{\sharp}_q = \{r\in \N: r \geq q\}
\]
satisfies 
\[
A = \bigcap_{q=1}^{\infty} A^{\sharp}_q = \emptyset 
\]
and
\[
hA = \bigcap_{q=1}^{\infty} hA^{\sharp}_q = \emptyset 
\]
 for all $h \in \N$ and so $H(A^{\sharp}_q)= \N$. 

If $S = \Z$ and $Q= \N$, then the sequence of sets 
$(A^{\sharp}_q)_{q=1}^{\infty}$ defined by 
\[
A^{\sharp}_q = \{r\in \Z: |r| \geq q\}
\]
satisfies 
\[
A = \bigcap_{q=1}^{\infty} A^{\sharp}_q = \emptyset. 
\]
By Theorem~\ref{intersect:theorem:Integers}, 
\[
hA = \emptyset \neq \Z =  \bigcap_{q=1}^{\infty} hA^{\sharp}_q 
\]
 for all $h \geq 2$ and so $H(A^{\sharp}_q)=\{1\}$.

Sum and product intersection sets other than \N\ and $\{1\}$ have been constructed.

\bt[Marques and Nathanson~\cite{marq-nath26}]       \label{intersect:theorem:mn-1}
In the additive semigroup \Z, for every positive integer $h_0$, 
there exists $(A_q)_{q=1}^{\infty} \in \mcf_{\N}^*(\Z)$ such that 
\[
H(A_q) = \{1\} \cup \{ h \in \N: h \geq h_0 \}.
\]
For every  integer $d \geq 2$,  there exists $(A_q)_{q=1}^{\infty} \in \mcf_{\N}^*(\Z)$ 
such that 
\[
H(A_q) = \{h \in \N: h \not\equiv 0 \pmod{d}\}.
\]
\et 

Very little is known about the arithmetic structure of product and sum intersection sets. 

We often use the following simple relation. 

\bl             \label{intersect:lemma:simple}
Let $S$ be a set, let $A$ and $B$ be subsets of $S$, 
and let $(B_q)_{q\in Q}$ be a family of subsets of $S$.   If 
\[
 \bigcap_{q\in Q} B_q = B
\]
then 
\[
 \bigcap_{q\in Q} (A \cup B_q) = A \cup B 
\]
and 
\[
 \bigcap_{q\in Q} (A \cap B_q) = A \cap B.
\]
\el

\begin{proof}
Klar.
\end{proof}

\section{Representation functions and intersection sets}

Let $S$ be a semigroup and let $A$ be a subset of  $S$.  
For all $x \in S$, we define the 
\emph{$h$-fold representation function} 
\[
r_{A,h}(x) = \card\left\{ (a_1,\ldots, a_h) \in A^h : a_1\cdots a_h = x \right\}.
\]
If $A \subseteq B$, then $r_{A,h}(x) \leq r_{B,h}(x)$ for all $x \in S$. 
We have 
\[
r_{A,1}(x) = 
\begin{cases}
1 & \text{if $x \in A$} \\
0 & \text{if $x \notin A$}
\end{cases}
\]
and $r_{A,h}(x) \geq 1$ if and only if $x \in A^h$.

Let $\Z^d$ be the additive group of vectors in $\R^d$ with integral coordinates, 
that is, the group of $d$-dimensional lattice points. 
Let $\N_0 = \{0,1,2,\ldots\}$ be the set of nonnegative integers 
and $\N_0^d$ the additive semigroup of $d$-dimensional lattice points with nonnegative integral coordinates.

\bl                \label{intersect:lemma:finiteness}
\benum
\item
In the multiplicative semigroup $\Z^* = \Z\setminus \{0\}$, 
\[
r_{\Z^*,h}(x) < \infty  \qquad \text{ for all  $h \geq 1$ and $x \in \Z^*$. }
\] 
\item 
In the additive semigroup $\N_0^d$, 
\[
r_{\N^d_0,h}(x) < \infty \qquad \text{ for all $h \geq 1$ and $x \in \N_0^d$.}
\]
\item
In the additive semigroup $\Z^d$, 
\[
r_{\Z^d,h}(x) = \infty \qquad \text{ for all $h \geq 2$ and $x \in \Z^d$.}
\]  
\eenum 
\el

\begin{proof}
(1) If $x$ is a nonzero integer and  and if $a_1,\ldots, a_h$ are nonzero integers 
such that $x = a_1\cdots a_h$, then each integer $a_i$ is a divisor of $x$ and 
every nonzero integer $x$ has only  finitely many positive or negative divisors.  
It follows that $r_{\Z^*,h}(x) < \infty$. 

(2) If $x = (u_1,\ldots, u_d)$ is a nonnegative lattice point and if $a_1,\ldots, a_h$ are  lattice points 
such that $x = a_1 +\cdots + a_h$, then each summand $a_i$ is in the finite hypercube 
\[
\left\{ (t_1,\ldots, t_d) \in \N_0^d : 0 \leq t_i \leq u_i \text{ for all } i = 1,\ldots, d \right\} 
\]
and so $r_{\N_0,h}(x) < \infty $. 

(3) For all lattice points $x$ and $y$ in $\Z^d$, we have 
\[
x = y + (x-y) +  \underbrace{0 + \cdots + 0}_{\text{$h-2$ summands}}
\] 
and so $r_{\Z^d,h}(x) = \infty $. 
This completes the proof. 
\end{proof}

The following result establishes a link between representation functions in number theory 
and intersection sets in semigroups.

\bt                              \label{intersect:theorem:finiteness}
Let $S$ be a semigroup and let $A$ be a subset of $S$.   
Let $(A_q)_{q=1}^{\infty}$ be a decreasing sequence of subsets of $S$ such that  
\[
A = \bigcap_{q=1}^{\infty} A_q. 
\] 
If $h \in \N$ and $r_{A_1,h}(x) < \infty$ for all $x \in A_1$, then $h \in H(A_q)$. \\
If $r_{A_1,h}(x) < \infty$ for all $h \in \N$ and $x \in A_1$, then $H(A_q) = \N$.
\et

\begin{proof} 
For all $h \in \N$, we have $A^h \subseteq \bigcap_{q=1}^{\infty} A_q^h$. 
Suppose that, for some $h \in \N$, we have $r_{A_1,h}(x) < \infty$ for all $x \in A_1$. 
We must prove that $ \bigcap_{q=1}^{\infty} A_q^h \subseteq  A^h$. 

If $x \in \bigcap_{q=1}^{\infty} A_q^h$, then, for all $q \in \N$, we have $x \in A_q^h$ 
and so there is an $h$-tuple $(a_{1,q},\ldots, a_{h,q})$ with $a_{i,q} \in A_q  \subseteq A_1$ 
for all $i = 1, \ldots, h$ 
such that 
\[
x = a_{1,q} \cdots a_{h,q}. 
\]
Because $r_{A_1,h}(x) < \infty$, there is an $h$-tuple $(b_1,\ldots, b_h )$ 
 with $b_i \in A_1$ for all $i = 1, \ldots, h$ and there is an infinite strictly increasing 
 sequence of positive integers $(q_j)_{j=1}^{\infty}$ such that 
\[
(a_{1,q_j}, \ldots, a_{h,q_j}) = (b_1,\ldots, b_h )
\]
for all $j = 1,2, 3, \ldots$.  
For all $q \in Q$, there exists $q_j$ such that $q \leq q_j$. 
Because the sequence  of sets $(A_q)_{q=1}^{\infty}$ is decreasing, we have 
$ A_{q_j} \subseteq A_q$ and so 
\[ 
(b_1,\ldots, b_h ) = (a_{1,q_j}, \ldots, a_{h,q_j}) \in A_{q_j} \times \cdots \times A_{q_j}
  \subseteq A_q \times \cdots \times A_q.  
 \] 
It follows that, for all $i = 1,\ldots, h$ and $q \in Q$, we have $b_i \in A_q$  and so 
$b_i \in \bigcap_{q=1}^{\infty} A_q = A$.  Therefore, $(b_1,\ldots, b_h ) \in A \times \cdots \times A$ 
and $x = b_1 \cdots b_h \in A^h$.  
This completes the proof.
\end{proof}

\bt                                     \label{intersect:theorem:finiteness-H}
If $S$ is a semigroup such that $r_{S,h}(x) < \infty$ for all $h \in \N$ and $x \in S$, then 
\[
\mathbf{H}^*_{\N}(S)  = \{ \N \}. 
\] 
In particular, for the multiplicative semigroups \N\ and $\Z^*$ 
and for the additive semigroup $\N_0^d$, 
\[
\mathbf{H}_{\N}(\N) = \mathbf{H}_{\N}(\Z^*) = \mathbf{H}_{\N}(\N_0^d) =\{ \N \}. 
\]
\et 

\begin{proof}
This follows from Lemma~\ref{intersect:lemma:finiteness}
and Theorem~\ref{intersect:theorem:finiteness}. 
\end{proof}

\bprob
Determine the semigroups $S$ such that 
\benum
\item 
for all $h \geq h_0$, there exists $x \in S$ with  $r_{S,h}(x) = \infty$, and  
\item
$\mbh_{\N}^*(S) = \{\N\}$.
\eenum
\eprob

The subset $A$ of a semigroup $S$ is a \emph{basis of order $h$} for $S$ 
if $A^h = S$ and a \emph{basis of exact order $h_0$} for $S$ if $h_0$ is the smallest 
integer such that $A^{h_0} = S$.  
 The subset $A$ is a \emph{nonbasis of order $h$} for $S$ 
if $A^h \neq S$ and a \emph{nonbasis} for $S$ if $A^h \neq S$ for all $h = 1, 2, 3, \ldots$. 

Let $S$ be a semigroup that contains an identity.
If $A$ is a basis of order $h_0$ for $S$, then $A$ is a basis of order $h$ for $S$ 
for all $h \geq h_0$.  If $A$ is a nonbasis of order $h_0$ for $S$, then $A$ is a nonbasis 
of order $h$ for $S$ for all $h \leq h_0$.

Additive bases and nonbases are described similarly.  
For example, the set $A = \{0\} \cup \{ \pm 2^i:i \in \N_0\}$ is a nonbasis 
for the additive abelian semigroup \Z. 

By Theorem~\ref{intersect:theorem:finiteness-H}, if $(A_q)_{q=1}^{\infty}$ is any  
asymptotically strictly decreasing sequence of sets of nonnegative integers, then 
$H(A_q) = \N$, that is, 
$hA = \bigcap_{q=1}^{\infty} hA_q$ for all $h \in \N$. 
This is not true for sets that contain both  positive and negative integers 
nor  for sets of nonnegative rational numbers.  
The following results  construct examples of such sets. 

\bt                       \label{intersect:theorem:Integers}
In the additive group \Z, for all $q \in \N$, let 
\[
 A^{\sharp}_q =  \{r\in \Z: |r| \geq q \}.  
\]
The sequence $(A^{\sharp}_q)_{q=1}^{\infty}$ is strictly decreasing and 
\[
 \bigcap_{q=1}^{\infty} A^{\sharp}_q = \emptyset.
\]
For all $h \geq 2$, 
\[
 \bigcap_{q=1}^{\infty} hA^{\sharp}_q = \Z
\]
and $H\left(A^{\sharp}_q\right) = \{1\}$.

Let $A$ be a co-infinite set of integers.  For all $q \in \N$, let 
\[
A_q =  A \cup A_q^{\sharp}. 
\]
The sequence $(A_q)_{q=1}^{\infty}$ is asymptotically strictly decreasing and 
\[
\bigcap_{q=1}^{\infty} A_q = A. 
\]
If $A$ is a basis of order $h$ for \Z, then $h \in H(A_q)$.  
If $A$ is a basis of exact order $h_0$ for \Z, then 
\[
H(A_q) = \{1\} \cup \{h \in \N_0:h \geq h_0 \}. 
\] 

If $A$ is a nonbasis of order $h_0$ for \Z, then 
\[
h \notin H(A_q) \qquad \text{ for all $h = 2,3,\ldots, h_0$.}
\]
If $A$ is a nonbasis for \Z, then $H(A_q) = \{1\}$.
\et

\begin{proof}
The sequence of sets $\left( A^{\sharp}_q  \right)_{q=1}^{\infty}$ is strictly decreasing and  
\[
 A^{\sharp} = \bigcap_{q=1}^{\infty} A^{\sharp}_q = \emptyset. 
\]
It follows that 
\beq                \label{intersect:2a}
h A^{\sharp} =  \emptyset  
\eeq
for all $h \in \N$. 

Let $h \geq 2$.   For  all $n \in \N_0$ and $q \in \N$, we have  
\[
\pm(n+q) \in A^{\sharp}_q 
\]
 and 
 \[
 \pm(h-1)q\in A^{\sharp}_q.  
 \] 
It follows that   
 \[
n =  (n+q) + \underbrace{q+ \cdots + q}_{\text{$h-2$ summands}}  
+ (-(h-1)q) \in h A^{\sharp}_q   
 \]
and 
\[
-n = -(n+q) +  \underbrace{ (-q) + \cdots +(-q)}_{\text{$h-2$ summands}} 
+ ((h-1)q)  \in h A^{\sharp}_q    
\]
and so $h A^{\sharp}_q = \Z$ for all  $q \in \N$.  Therefore, 
\beq                \label{intersect:2b}
\bigcap_{q=1}^{\infty} h A^{\sharp}_q = \Z. 
\eeq 
From~\eqref{intersect:2a} and~\eqref{intersect:2b}, we obtain 
\[
h A^{\sharp} \neq \bigcap_{q=1}^{\infty} h A^{\sharp}_q 
\]
for all $h \geq 2$ and so $H\left(A_q^{\sharp} \right) = \{ 1\}$.

Because the set $A_q$ is co-finite and the set $A$ is co-infinite, the set $A_q\setminus A $ 
is infinite for all $q \in \N$.  It follows that, for all $q \in Q$, 
there exists 
\[
s \in A_q\setminus A \subseteq A_q.
\] 
We have  $|s| \geq q$ and so $|s| +1 > q$.  Therefore,  
\[
s \notin A_{|s| +1}, \qquad  s \in A_q, \qqand  A_{|s| +1} \neq  A_q 
\]
and so the sequence $(A_q)_{q=1}^{\infty}$ is asymptotically strictly decreasing. 

For all $q \in \N$, we have $A \subseteq A_q$ and so 
$A \subseteq \bigcap_{q=1}^{\infty} A_q$. 
Let $x \notin A$.  For all $q > |x|$, we have $x \notin A^{\sharp}_q$ and so $x \notin A \cup A^{\sharp}_q = A_q$. 
Therefore, 
\[
x \notin \bigcap_{q=1}^{\infty} A_q  
\]
and 
\[
A =  \bigcap_{q=1}^{\infty} A_q.
\]
For all $h \geq 2$, we have 
$\Z = hA_q^{\sharp} \subseteq hA_q \subseteq \Z$ and so 
\[
\bigcap_{q=1}^{\infty} hA_q = \Z. 
\]  If $A$ is a basis of  order $h $ for \Z, then 
$hA = \Z = \bigcap_{q=1}^{\infty} hA_q$ and $h \in H(A_q)$.  
If $A$ is a basis of exact order $h_0$, then 
$hA = \Z = \bigcap_{q=1}^{\infty} hA_q$ if and only if $h \geq h_0$ and  
\[
H(A_q) = \{1\} \cup \{h \in \N: h \geq h_0 \}. 
\]

If $A$ is a nonbasis of order $h_0$ for \Z, then, for all $h = 2,\ldots, h_0$, 
we have $hA \neq \Z =  \bigcap_{q=1}^{\infty} hA_q$   
and $h \notin H(A_q)$. 

If $A$ is not a basis  for \Z, then $hA \neq \Z = \bigcap_{q=1}^{\infty} hA_q$ 
for all $h \geq 2$ and so 
\[
H(A_q) = \{1\}. 
\] 
This completes the proof. 
\end{proof}

\bt         \label{intersect:theorem:rational} 
There exist asymptotically  strictly decreasing sequences $(A_q)_{q=1}^{\infty}$ of sets of 
positive rational numbers with $A = \bigcap_{q=1}^{\infty} A_q$ and   
\[
hA \neq \bigcap_{q=1}^{\infty} hA_q
\]
for all $h \geq 2$, that is, $H(A_q) = \{1\}$. 
\et

\begin{proof}
Let $B^* = \{b_n: n \in \N\}$ be a set of positive integers such that 
\[
b_1 > 1 \qqand b_{n+1} > b_n + 2
\]
for all $n \in \N$.  
For all $q \in \N$, let 
\[
B_{n,q} =  \left\{ b_n + \frac{1}{r}: |r| \geq q\right\}. 
\]
The inequality  
\[
\frac{2}{q} \leq 2 <b_{n+1} - b_n 
\]
implies 
\[
\max B_{n,q} = b_n + \frac{1}{q} < b_{n+1}-\frac{1}{q} = \min B_{n+1,q} 
\]
and so $(B_{n,q})_{n=1}^{\infty}$ is a sequence of pairwise disjoint 
sets of positive rational numbers.  
Let  
\[
A_q = \bigcup_{n=1}^{\infty} B_{n,q}. 
\]
For all $n \in \N$, we have 
\[
\bigcap_{q=1}^{\infty} B_{n,q} = \emptyset  
\]
and so $(A_q)_{q=1}^{\infty}$ is a strictly decreasing sequence of sets 
of positive rational numbers such that  
\[
A = \bigcap_{q=1}^{\infty} A_q = \emptyset.  
\] 
It follows that 
\beq          \label{intersect:hB-0} 
hA = \emptyset 
\eeq
for all $h \geq 2$.  
We shall prove that 
\beq          \label{intersect:hB} 
 \bigcap_{q=1}^{\infty} hA_q  = hB^* 
\eeq
for all $h \geq 2$.

The identity  
\[
\frac{1}{q} - \frac{1}{q} = \frac{1}{q} - \frac{1}{2q} - \frac{1}{2q} = 0
\]
 implies that, for all $h \geq 2$,  there exist $\varepsilon_i \in \{1,-1\}$ 
and $\nu_i \in \{1,2\}$ such that 
\[
\sum_{i=1}^h \frac{\varepsilon_i}{\nu_i q} = 0.
\]

If $x \in hB^*$, then there exist $b_{n_1},\ldots, b_{n_h} \in B^*$ 
such that 
\[
x = b_{n_1}  + \cdots + b_{n_h}. 
\]
We have 
\[
b_{n_i}+  \frac{\varepsilon_i}{\nu_i q} \in A_q 
\]
for all $i = 1,\ldots, h$  and so 
\[
x = b_{n_1}  + \cdots + b_{n_h} = \left( b_{n_1} +  \frac{\varepsilon_1}{\nu_1q}  \right)  
+ \cdots + \left( b_{n_h} +  \frac{\varepsilon_h}{\nu_h q}  \right)  \in hA_q
\]
for all $q \in \N$.  Therefore, $x \in \bigcap_{q=1}^{\infty} hA_q$ and 
\beq                                       \label{intersect:6a} 
hB^* \subseteq  \bigcap_{q=1}^{\infty} hA_q. 
\eeq

Let $h \geq 2$ and $x \in \bigcap_{q=1}^{\infty} hA_q$.  
Choose an  integer $q$ such that 
\[
2h < q. 
\]
Because $x \in hA_q$,  there exist integers $b_{n_i}$ and  $r_i$ with $|r_i| \geq q$  and 
\[
b_{n_i} +  \frac{1}{r_i} \in A_q 
\]
such that 
\[
x = \sum_{i=1}^h \left( b_{n_i} + \frac{1}{r_i} \right) 
= \sum_{i=1}^h b_{n_i} + \sum_{i=1}^h  \frac{1}{r_i}. 
\]
For all $q ' > q$, we also have $x \in hA_{q'}$ 
 and  there   exist integers $b_{m_j}$ and $r'_j$ with $|r'_j| \geq q'$ and 
\[
b_{m_j} +  \frac{1}{r'_j} \in A_{q'} 
\]
such that 
\beq             \label{intersect:xy}
x = \sum_{j=1}^h \left( b_{m_j} + \frac{1}{r'_j} \right) 
=  \sum_{j=1}^h b_{m_j} + \sum_{j=1}^h  \frac{1}{r'_j}. 
\eeq
Then 
\begin{align*}
\left| \sum_{i=1}^h   b_{n_i} -  \sum_{j=1}^h b_{m_j} \right| 
& =  \left|   \sum_{j=1}^h  \frac{1}{r'_j}  -  \sum_{i=1}^h  \frac{1}{r_i} \right|  \\
& \leq   \sum_{j=1}^h  \frac{1}{ \left|   r'_j  \right|}  +  \sum_{i=1}^h  \frac{1}{ \left|   r_i\right| } \\
& \leq \frac{h}{q'} + \frac{h}{q} < \frac{2h}{q} < 1. 
\end{align*}

Because  $ \sum_{i=1}^h b_{n_i} $ and $ \sum_{j=1}^h b_{m_j} $ are integers, we have 
\[
\sum_{i=1}^h   b_{n_i} = \sum_{j=1}^h b_{m_j} = y \in hB^*. 
\]
It follows from~\eqref{intersect:xy} that 
\[
|x-y| = \left|   \sum_{j=1}^h  \frac{1}{r'_j} \right| \leq \frac{h}{q'}
\]
for all $q' > q$, and so  $x = y \in hB^*$.  Therefore, 
\beq                                       \label{intersect:6b} 
\bigcap_{q=1}^{\infty} hA_q \subseteq hB^*  
\eeq
for all $h \geq 2$.    Relations~\eqref{intersect:6a} 
and~\eqref{intersect:6b} imply~\eqref{intersect:hB}  and so $H(A_q) = \{1\}$. 
This completes the proof. 
\end{proof}

In the proof of Theorem~\ref{intersect:theorem:rational}, if we define 
\[
B'_{n,q} = \{b_n\} \cup  \left\{ b_n + \frac{1}{r}: |r| \geq q\right\} 
\]
and 
\[
A'_q = \bigcup_{n=1}^{\infty} B'_{n,q} = A_q \cup B^*  
\]
then we obtain a strictly decreasing sequence $(A'_q)_{q=1}^{\infty}$ 
of sets of positive rational numbers such that $\bigcap_{q=1}^{\infty} hA'_q = hB^*$ 
for all $h \in \N$ and so $H(A'_q) = \N$.

\bt         \label{intersect:theorem:open} 
There exist asymptotically  strictly decreasing sequences $(A_q)_{q=1}^{\infty}$ of sets 
of real numbers with $A = \bigcap_{q=1}^{\infty} A_q$ such that each set $A_q$ 
is a union of open intervals and 
\[
hA \neq \bigcap_{q=1}^{\infty} hA_q
\]
for all $h \geq 2$, that is, $H(A_q) = \{1\}$. 
\et

\begin{proof}
Let $B^* = \{b_n: n \in \N \}$ be a set of real numbers such that 
such that 
\[
b_{n+1} > b_n + 2
\]
for all $n \in \N$.  
For all $q \in \N$, let 
\[
B_{n,q} =  \left( b_n - \frac{1}{q}, b_n \right) 
\cup  \left( b_n, b_n + \frac{1}{q} \right)
\]
As in the proof of Theorem~\ref{intersect:theorem:rational},
the inequality  
\[
\frac{2}{q} \leq 2 <b_{n+1} - b_n 
\]
implies 
\[
\max B_{n,q} = b_n + \frac{1}{q} < b_{n+1}-\frac{1}{q} = \min B_{n+1,q} 
\]
and so $(B_{n,q})_{n=1}^{\infty}$ is a sequence of pairwise disjoint sets, each a union 
of two open intervals.   Let  
\[
A_q = \bigcup_{n=1}^{\infty} B_{n,q}. 
\]
For all $n \in \N$, we have 
\[
\bigcap_{q=1}^{\infty} B_{n,q} = \emptyset  
\]
and so $(A_q)_{q=1}^{\infty}$ is a strictly decreasing sequence of sets  such that  
\[
A = \bigcap_{q=1}^{\infty} A_q = \emptyset.  
\] 
It follows that 
\[
hA = \emptyset 
\]
for all $h \geq 2$.  
The proof that 
\[
 \bigcap_{q=1}^{\infty} hA_q  = hB^* 
\]
for all $h \geq 2$ is exactly the same as the proof in Theorem~\ref{intersect:theorem:rational}.
\end{proof}

In the proof of Theorem~\ref{intersect:theorem:open}, if we define 
\[
B'_{n,q} =  \left( b_n - \frac{1}{q},  b_n + \frac{1}{q} \right)
\]
and 
\[
A'_q = \bigcup_{n=1}^{\infty} B'_{n,q} = A_q \cup B^*  
\]
then we obtain a strictly decreasing sequence $(A'_q)_{q=1}^{\infty}$ 
of open sets of real numbers such that $\bigcap_{q=1}^{\infty} hA'_q = hB^*$ 
for all $h \in \N$ and so $H(A'_q) = \N$. 
The sets $B_{n,q}$ and $B'_{n,q}$ are sets of positive real numbers if we choose $b_1 > 1$.

\section{Product intersection sets in groups} 
Every infinite abelian group contains an infinite proper subgroup.  
This is not necessarily true for infinite nonabelian groups.  
Of course, if an infinite group $G$ is not a torsion group, then $G$ contains an element $g$ 
of infinite order and the cyclic subgroup generated by $g^2$ is an infinite 
proper subgroup of $G$.  
Infinite torsion groups, however, do not always contain infinite proper subgroups. 
For example, for every prime number $p$, there exist infinite groups, 
called ``Tarski monsters'', such that every proper subgroup is finite and 
cyclic of order $p$.  

\bt                        \label{intersect:theorem:subgroup} 
Let $G$ be an infinite group, not necessarily abelian, that contains an infinite proper subgroup.  
There exists a  strictly decreasing sequence $(A_q)_{q=1}^{\infty}$ 
of subsets of $G$ such that $H(A_q) = \{1\}$. 
\et 

\begin{proof} 
Let $A$ be an infinite proper subgroup of $G$.  Then 
\[
A^h = A
\]
for all $h \geq 1$. 
Choose  $x \in G\setminus A$.  Because $A$ is infinite, 
the coset $xA$ is infinite and disjoint from $A$.  
Let $(x_q)_{q=1}^{\infty}$ be an infinite sequence of distinct elements of $xA$.  
Then 
\[
x_qA = xA
\]
for all $q \in N$. 
We define the sets 
\[
A_q^{\sharp} = \{x_r: r \geq q\}
\] 
and 
\[
A_q = A \cup A_q^{\sharp}. 
\]
Then 
\[
\bigcap_{q=1}^{\infty} A_q^{\sharp} = \emptyset 
\qqand 
\bigcap_{q=1}^{\infty} A_q= A. 
\]
For all $h \geq 2$ and $q \in Q$, we have  $A \cup \{x_q\} \subseteq A_q$  
and so 
\[
A_q^h \supseteq (A \cup \{x_q\})^h \supseteq A^h \cup x_q A^{h-1} 
 = A \cup x_q A = A \cup x A. 
\] 
Therefore,  
\[
A^h = A \neq A  \cup xA  \subseteq \bigcap_{q=1}^{\infty} A_q^h    
\]
and so $h \notin H(A_q)$ for all $h \geq 2$.  Thus, $H(A_q) = \{ 1 \}$. 
This completes the proof. 
\end{proof}

\bl              \label{intersect:lemma:surjection}
Let $G$ and $G'$ be groups, not necessarily abelian and let  
 $B$ be a nonempty subset of $G'$.   
If $f:G\rightarrow G'$ is a surjective group homomorphism, then 
\[
 \left( f^{-1} (B) \right)^h =  f^{-1} \left(B^h\right) 
\]
for all $h \in \N$. 
\el

\begin{proof}
Let  $A = f^{-1}(B)$. 
If $x \in A^h = \left( f^{-1} (B) \right)^h$, then there exist 
$a_1,\ldots, a_h \in A $ such that 
\[
x = a_1\cdots a_h
\]
and so 
\[
f(x) = f(a_1)\cdots f(a_h) \in B^h
\]
and $x \in f^{-1}\left( B^h\right)$.  Thus,
\[
 \left( f^{-1} (B) \right)^h \subseteq  f^{-1} \left(B^h\right). 
\]

Conversely, if $x \in f^{-1} \left(B^h\right)$, then $f(x) \in B^h$ and there exist 
$b_1,\ldots, b_h \in B$ such that 
\[
f(x) = b_1\cdots b_h.
\]  
Because $f$ is surjective, there exist $a_i \in A$ such that $f(a_i) = b_i$ 
for all $i = 1,\ldots, h-1$.  Let
\[
a_h = a_{h-1}^{-1} \cdots a_2^{-1}a_1^{-1}x \in G.
\] 
Then 
\begin{align*}
f(a_h) & =  f(a_{h-1})^{-1} \cdots  f(a_1)^{-1}  f(x) 
= b_{h-1}^{-1} \cdots b_1^{-1}  f(x) = b_h \in B 
\end{align*} 
and so $a_h \in f^{-1}(B) = A$ and 
\[
x = a_1\cdots a_h \in A^h = \left(f^{-1}(B)\right)^h. 
\]
Therefore, 
\[
f^{-1} \left(B^h\right) \subseteq  \left( f^{-1} (B) \right)^h.  
\] 
This completes the proof. 
\end{proof}

\bt
Let $G$ and $G'$ be groups.  
If there exists a surjective group homomorphism $f:G\rightarrow G'$, 
then  $\mathbf{H}^*_{\N} (G') \subseteq \mathbf{H}^*_{\N}(G)$ 
and $\mathbf{H}_Q (G') \subseteq \mathbf{H}_Q(G)$ for all index sets $Q$. 
\et

\begin{proof}
Let $(B_q)_{q\in Q}$ be a family of subsets of $G'$ 
and let $(A_q)_{q\in Q}$ be the family of subsets of $G$ 
defined by $A_q = f^{-1}(B_q)$ for all $q \in Q$.  
Because $f$ is surjective, we have $f(A_q) = B_q$ for all $q \in Q$. 

Let 
\[
A = \bigcap_{q\in Q} A_q  \qqand B = \bigcap_{q\in Q} B_q. 
\]
Then 
\[
A = \bigcap_{q\in Q} A_q = \bigcap_{q\in Q}  f^{-1}(B_q) 
= f^{-1} \left(  \bigcap_{q\in Q}   B_q\right) = f^{-1} \left( B \right) 
\] 
and, because $f$ is surjective, we have  $f(A) = B$. 

We apply Lemma~\ref{intersect:lemma:surjection}. 
If $h \in H(B_q)$, then $B^h = \bigcap_{q\in Q} B_q^h$ and 
\begin{align*}
A^h & = \left(f^{-1}(B)\right)^h =  f^{-1}(B^h) \\
& = f^{-1}\left( \bigcap_{q\in Q} B_q^h\right) =  \bigcap_{q\in Q} f^{-1}\left(  B_q^h\right) \\ 
& =  \bigcap_{q\in Q} \left( f^{-1}\left(  B_q\right) \right)^h 
=  \bigcap_{q\in Q} A_q^h 
\end{align*}
and so $h \in H(A_q)$.  This proves $H(B_q) \subseteq H(A_q)$. 

If $h \in H(A_q)$, then $A^h = \bigcap_{q\in Q} A_q^h$ and 
\begin{align*}
f^{-1}\left( B^h \right) & =  \left( f^{-1}\left( B \right) \right)^h = A^h \\
& =  \bigcap_{q\in Q} A_q^h = \bigcap_{q\in Q}  \left( f^{-1}\left( B_q \right) \right)^h \\
& = \bigcap_{q\in Q}    f^{-1}  \left( B_q ^h \right)  =     f^{-1}  \left(  \bigcap_{q\in Q}   B_q ^h \right). 
\end{align*} 
Because $f$ is surjective, we have  
\[
B^h = f\left( f^{-1}\left( B^h \right) \right) =  f\left(   f^{-1}  \left(  \bigcap_{q\in Q}   B_q ^h \right) \right) =  \bigcap_{q\in Q}   B_q ^h 
\]
and so  $h \in H(B_q)$ and $H(A_q) \subseteq H(B_q)$.    
This proves $H(A_q) = H(B_q)$ and so 
$\mathbf{H}_Q (G') \subseteq \mathbf{H}_Q(G)$.

Because $f$ is surjective, if the sequence $(B_q)_{q=1}^{\infty} $ is asymptotically 
strictly decreasing, then the sequence $ \left( f^{-1} \left( B_q \right) \right)_{q=1}^{\infty}$ 
is asymptotically strictly decreasing and 
$\mathbf{H}^*_{\N} (G') \subseteq \mathbf{H}^*_{\N}(G)$. 
This completes the proof. 
\end{proof}

\bl         \label{intersect:lemma:co-finite-basis} 
Let $G$ be an infinite group.  
Every co-finite subset $X$ of $G$ is a basis of order $h$ for $G$ for all $h \geq 2$.
\el 

\begin{proof}
We must show that $g \in X^h$ for all  $g \in G$ and $h \geq 2$.  

If $X$ is a co-finite subset of $G$, then the set $ X^{-1} g = \{ y^{-1} g :y \in X\}$ is co-finite.  
Let $x_3,\ldots, x_h \in X$.  The set $Xx_3\cdots x_h = \{zx_3\cdots x_h: z \in X\}$ is also co-finite. 
Because the intersection of co-finite subsets of an infinite set is co-finite, we have 
\[
X^{-1} g  \cap Xx_3\cdots x_h  \neq \emptyset 
\]
and so there exist $x_1,x_2 \in X$ such that $x_1^{-1} g = x_2x_3\cdots x_h$ 
and $g = x_1x_2\cdots x_h \in X^h$. 
This completes the proof. 
\end{proof}

\bt             \label{intersect:theorem: sharp}
Let $G$ be an infinite group and let $(A_q^{\sharp})_{q\in Q}$ be a family 
of co-finite subsets of $G$.  Let $A$ be a subset of $G$ and let 
\[
A_q = A \cup A_q^{\sharp} 
\]
for all $q \in Q$.  Then 
\[
H(A_q) = \{1\} \cup \{h \in \N: A^h = G\}.
\]
\et

\begin{proof} 
 If $(A_q^{\sharp})_{q\in Q}$ is a family of co-finite subsets of $G$, 
 then, by Lemma~\ref{intersect:lemma:co-finite-basis}, 
 $\left( A_q^{\sharp}\right)^h = G$ for all $q \in Q$ and so 
 $\bigcap_{q\in Q} \left( A_q^{\sharp}\right)^h = G$. 
 Because $A_q^{\sharp} \subseteq A_q$ for all $q \in Q$, we have 
  $\bigcap_{q\in Q} A_q^h = G$ and so 
  $h \in H(A_q)$ if and only if $h=1$ or $A^h = G$.
This  completes the proof. 
 \end{proof}

\bt            \label{intersect:theorem: sharp-strict}
Let $G$ be an infinite group and let $(A_q^{\sharp})_{q=1}^{\infty}$ be an asymptotically 
strictly decreasing sequence of co-finite subsets of $G$ such that 
\[
\bigcap_{q=1}^{\infty}A_q^{\sharp} = \emptyset. 
\] 
Let $A$ be a co-infinite subset of $G$ and let 
\[
A_q = A \cup A_q^{\sharp} 
\]
for all $q \in Q$.  Then $(A_q)_{q=1}^{\infty}$ is an asymptotically strictly 
decreasing sequence of co-finite subsets of $G$ and 
\[
H(A_q) = \{1\} \cup \{h \in \N: A^h = G\}.
\]
\et

\begin{proof} 
The sequence $(A_q)_{q=1}^{\infty}$ is decreasing because the sequence
$(A_q^{\sharp})_{q=1}^{\infty}$ is decreasing.   
For all $q \in Q$, the set $A_q^{\sharp} \setminus A$ is infinite because $A^{\sharp}_q$ 
is a co-finite subset of the infinite group $G$ and $A$ is co-infinite,   
Because  $(A_q)_{q=1}^{\infty}$ is decreasing and $\bigcap_{q=1}^{\infty}A_q^{\sharp} = \emptyset$, 
for all $x \in A_q^{\sharp} \setminus A$, there exists $q' > q$ such that 
$x \notin A_{q'}^{\sharp}$ and so $x \notin A_{q'}^{\sharp} \setminus A$.  
Thus,  $A_{q'} = A \cup A_{q'}^{\sharp}$ is a proper subset of $A_q = A \cup A_q^{\sharp}$ 
and the sequence $(A_q)_{q=1}^{\infty}$ is asymptotically strictly decreasing.  
The relation $H(A_q) = \{1\} \cup \{h \in \N: A^h = G\}$ 
follows from Theorem~\ref{intersect:theorem: sharp}.
\end{proof}

\bt
Let $A$ be a nonempty finite set of integers.  There is an asymptotically strictly decreasing 
sequence $(A_q)_{q=1}^{\infty}$ of sets of integers such that 
\[
hA = \bigcap_{q=1}^{\infty} hA_q
\] 
for all $h \in \N$.  Equivalently,  $\N\in \mbh_{\N}^*(A)$.
\et

\begin{proof} 
We begin  with the observation that  if $m_q$ and $m_{q'}$ are positive integers 
such that $m_q < m_{q'}$ and $m_q$ divides $m_{q'}$, then, for every integer $a$,  
the congruence class $a \pmod{m_{q'}}$ is 
a proper subset of the congruence class  $a \pmod{m_q}$. 

Because $A$ is finite, we can choose an integer $m^*$ such that $A \subseteq [-m^*,m^*]$.  
Let ($m_q)_{q=1}^{\infty}$ be a strictly increasing sequence of integers with 
$m_1 > 2m^*$ such that  $m_q$ divides $m_{q+1}$ for all $q \geq 1$.  
Let 
\[
A_q =  \left\{   n \in \Z: n \equiv a \pmod{m_q} \text{ for some } a \in A \right\} 
= \bigcup_{a\in A} a \pmod{m_q}.
\]
If $a,a' \in A$ and $a \equiv a' \pmod{m_q}$, then $m_q \geq m_1 > 2m^*$ implies $a = a'$.  
It follows that  $A_q$ is the union of $|A|$ pairwise disjoint congruence classes modulo $m_q$. 
Because 
\[
A_{q+1} =  \bigcup_{a\in A} a \pmod{m_{q+1}}
\]
and the congruence class $a \pmod{m_{q+1}}$ is a proper subset of the class $a \pmod{m_q}$ 
for all $a \in A$, 
it follows that $A_{q+1}$ is a proper subset of $A_q$ and so the sequence 
$(A_q)_{q=1}^{\infty}$ is strictly decreasing. 
 
For all $q,h \in \N$, we have $A \subseteq A_q$ and $hA \subseteq hA_q$ and so 
\[
A \subseteq \bigcap_{q=1}^{\infty} A_q 
\qqand 
hA \subseteq \bigcap_{q=1}^{\infty} hA_q. 
\]

Let $h \in \N$ and  $x \in \bigcap_{q=1}^{\infty} hA_q$.
Because $\lim_{q\rightarrow\infty} m_q = \infty$, we can choose $q$ such that $m_q > 2hm^*$.  
Then $hA \subseteq [-hm^*,hm^*]$ implies that the elements of the sumset $hA$ are 
pairwise incongruent modulo $m_q$. 
For all $q' > q$, we have $x \in hA_q$ and $x \in hA_{q'}$, 
and so there exist $b,b' \in hA$ such that 
\[
x\equiv b \pmod{m_q}
\]
and 
\[
x\equiv b' \pmod{m_{q'}}.
\]
Because $m_q$ divides $m_{q'}$, we have  
\[
x\equiv b' \pmod{m_q}
\]
and so 
\[
b \equiv b' \pmod{m_q}.  
\]
Because the elements of $hA$ are pairwise incongruent modulo $m_q$, it follows 
that $b=b'$ and so 
\[
x \equiv b \pmod{m_{q'}}
\]
for all $q' \geq q$.  Because $\lim_{q\rightarrow\infty} m_q = \infty$, 
we have $x = b \in hA$ and so $hA = \bigcap_{q=1}^{\infty} hA_q$.  
This completes the proof. 
\end{proof}

The usual Euclidean length of a vector $x_i = (x_{i,1},\ldots, x_{i,d})$ in $\R^d$ 
is $\|x_i \| = \sqrt{\sum_{j=1}^d x_{i,j}^2}$. 
The vector $x_i$ is \emph{nonnegative} if its coordinates are nonnegative, that is, 
$x_{i,j} \geq 0$ for all $j = 1,\ldots, d$. 
If $x_1,\ldots, x_k$ are nonzero vectors such that $\sum_{i=1}^k x_i = 0$, then 
\[
\|x_1  + \cdots + x_k \| <  \min(\|x_1\|,  \ldots, \|x_k\| ).
\]
This inequality is not possible for nonnegative vectors.

\bl            \label{intersect:lemma:vectorMin}
If  $x_1,  \ldots, x_k$ are nonnzero nonnegative vectors in $\R^d$, then 
\[
\|x_1  + \cdots + x_k \|   \geq \sqrt{k}  \min(\|x_1\|,  \ldots, \|x_k\| ).
\] 
\el

\begin{proof}
Let $x_i = (x_{i,1},\ldots, x_{i,d})$ for $i = 1,\ldots, k$.  
The coordinates $x_{i,j}$ are nonnegative, and so 
\begin{align*}
\| x_1+\cdots + x_k \|^2 & = \sum_{j=1}^d \left( \sum_{i=1}^k x_{i,j} \right)^2 \\
& \geq \sum_{j=1}^d  \sum_{i=1}^k x_{i,j}^2  = \sum_{i=1}^k  \sum_{j=1}^d x_{i,j}^2    \\ 
& = \sum_{i=1}^k \|x_i\|^2 \\
& \geq k   \min(\|x_1\|, \ldots, \|x_k\| )^2
\end{align*} 
and so 
\[
\|x_1+ \cdots + x_k \|  \geq \sqrt{k}  \min(\|x_1\|, \ldots, \|x_k\| ).
\]
This completes the proof. 
\end{proof}

\bt
Let $A$ be a nonempty finite set of lattice points in $\Z^d$.  
There is an asymptotically strictly decreasing 
sequence $(A_q)_{q=1}^{\infty}$ of sets of lattice points  such that 
\[
hA = \bigcap_{q=1}^{\infty} hA_q
\] 
for all $h \in \N$.  Equivalently,  $\N\in \mbh_{\N}^*(A)$.
\et

\begin{proof}
Let $A$ be a finite set of lattice points in $\Z^d$.  
The set $A$ is bounded.  
Choose $m^* > 0$ such that $\|a\| \leq m^*$ for all $a \in A$.  
For all $q \in \N$, let 
\[
A_q^{\sharp} = \{x \in \N_0^d: \| x\| \geq 2qm^*\} 
\]
and
\[
A_q  = A \cup A_q^{\sharp}  = A \cup \{x \in \N_0^d: \| x\| \geq 2qm^*\}. 
\]
The sequences $(A_q^{\sharp})_{q=1}^{\infty}$ and $(A_q)_{q=1}^{\infty}$ are 
strictly decreasing, and $A \cap A_q^{\sharp} = \emptyset$ for all $q \in \N$. 
We have 
\[
\bigcap_{q=1}^{\infty} A_q^{\sharp} = \emptyset  
\]
and so
\[
  \bigcap_{q=1}^{\infty}  A_q = A. 
\]
We shall prove that 
\beq                                                      \label{intersect:bigcap-all h}
hA = \bigcap_{q=1}^{\infty} hA_q
\eeq
for all $h \in \N$. 

We have $hA \subseteq  \bigcap_{q=1}^{\infty}  hA_q$.  
If  $hA \neq  \bigcap_{q=1}^{\infty}  hA_q$, then there is a lattice point 
 $x \in \bigcap_{q=1}^{\infty}  hA_q$ such that $x \notin hA$.  
It follows that, for all $q \in \N$,  there exist  
$k \in \{1,2,\ldots, h \}$ and lattice points 
$a_{1,q},\ldots, a_{k,q} \in A_q \setminus A = A_q^{\sharp} $  
and  lattice points  $a_1,\ldots, a_{h-k} \in A$ such that 
\[
x  = \sum_{i=1}^k a_{i,q} + \sum_{i=1}^{h-k} a_i.  
\]
The lattice points $a_{i,q}$ are nonnegative and nonzero.  
Applying  Lemma~\ref{intersect:lemma:vectorMin}, we obtain 
\begin{align*}
\left\| x - \sum_{i=1}^{h-k} a_i \right\| & =  \left\| \sum_{i=1}^k a_{i,q} \right\| \\
& > \min\left( \|  a_{1,q} \|, \ldots,  \|  a_{k,q} \| \right) \\
& \geq 2qm^*. 
\end{align*}
We also have 
\begin{align*}
\left\| x - \sum_{i=1}^{h-k} a_i \right\| & \leq \| x\| +   \sum_{i=1}^{h-k} \left\|  a_i\right\| \\
& < \| x\| + hm^*
\end{align*}
and so 
\[
2qm^* < \| x\| + hm^*
\]
for all $q \in \N$, which, because $m^* > 0$,  is absurd. 
This completes the proof.  
\end{proof}

\bt
Let $A$ be a nonempty finite set of a countably infinite group $G$. 
There is an asymptotically strictly decreasing 
sequence $(A_q)_{q=1}^{\infty}$ of subsets of $G$ with 
$A = \bigcap_{q=1}^{\infty} A_q$ such that $H(A_q)=1$. 
Equivalently,  $\{ 1\} \in \mch_{\N}^*(A)$.
\et

\begin{proof}
The set $G\setminus A$ is countably infinite.  Let 
\[
G \setminus A = \{a_q:q=1,2,\ldots\}
\]
where $a_q \neq a_r$ for all $r\neq q$. 
For all $q \in \N$, let 
\[
A_q^{\sharp} = \{a_r: r \geq q\}
\]
and 
\[
A_q = A \cup A_q^{\sharp}. 
\]
The sequences of sets $(A^{\sharp}_q)_{q=1}^{\infty}$  and $(A_q)_{q=1}^{\infty}$
are strictly decreasing with 
\[
\bigcap_{q=1}^{\infty} A^{\sharp}_q  = \emptyset 
\qqand \bigcap_{q=1}^{\infty} A_q  = A.
\]
The set $A_q$ is co-finite.  By Lemma~\ref{intersect:lemma:co-finite-basis},  $A_q^h = G$ for all $h \geq 2$.  
However, $A$ is finite and so $A^h$ is finite for all $h \in \N$. 
Thus, for all $h \geq 2$, we have 
\[
A^h \neq G = \bigcap_{q=1}^{\infty} A_q^h 
\]
and so 
\[
H(A_q) = \{1\} \in \mch_{\N}^*(A) \subseteq \mbh_{\N}^*(G).
\]
This completes the proof. 
\end{proof}

\bprob
Let $G$ be an infinite torsion group.  
Describe the product intersection sets $\mbh_Q(G)$. 
Must there exist $X \in \mbh_{\N}^*(G)$ with $X \neq \{1\}$ and $X \neq \N$? 
\eprob

\section{Intersections of product intersection sets} 
In this section we prove a global intersection relation for product intersection sets.

\bl           \label{intersect:lemma:times}
Let $A$ be a subset of a semigroup $S_1$ and $B$ a subset of a semigroup $S_2$.  
Then $A \times B$ is a subset of the product semigroup $S_1 \times S_2$ and, 
for all $h \in \N$, 
\[
A^h \times B^h = (A\times B)^h. 
\] 
\el

\begin{proof}
If $(x,y) \in A^h \times B^h$, then $x \in A^h$ and $y \in B^h$ and so there exist 
$a_1,\ldots, a_h \in A$ and $b_1,\ldots, b_h \in B$ such that 
$x = a_1\cdots a_h$ and $y = b_1\cdots b_h$.  
We have $(a_i,b_i) \in A \times B$ for $i = 1,\ldots, h$ and 
\[
(x,y) = (a_1,b_1) \cdots (a_h,b_h) \in (A\times B)^h.  
\] 
Thus, 
\[
 A^h \times B^h \subseteq (A\times B)^h.  
\]

Conversely, if $(x,y) \in  (A\times B)^h$, then there exist 
$(a_1,b_1), \ldots, (a_h,b_h) \in A\times B$ such that 
\[
(x,y) = (a_1,b_1) \cdots (a_h,b_h) = (a_1\cdots a_h, b_1\cdots b_h). 
\]
We have $x = a_1\cdots a_h\in A^h$ and $y = b_1\cdots b_h \in B^h$ 
and so $(x,y) \in A^h \times B^h$ and 
\[
(A\times B)^h \subseteq A^h \times B^h. 
\]
This completes the proof. 
\end{proof}

\bt
Let $\mathbf{Semi}$ be the class of all semigroups.  
For every index set $Q$, let 
\[
\mbh_Q = \bigcup_{S\in \mathbf{Semi}} \mbh_Q(S) 
=  \bigcup_{S\in \mathbf{Semi}} \  \bigcup_{(A_q) \in \mcf_Q(S)} H^S_Q(A_q)
\] 
be the set of all product intersection sets $H_Q(S)$ 
for all semigroups $S$.  Let 
\[
\mbh_{\N}^* = \bigcup_{S\in \mathbf{Semi}} \mbh_{\N}^*(S). 
\] 
  The sets $\mbh_Q$ and $\mbh_{\N}^*$ are closed under intersection. 
\et

\begin{proof} 
Let $S_1$ and $S_2$ be semigroups and let $S_1\times S_2$ be the product semigroup.  
For all $(A_q)_{q\in Q} \in \mcf_Q(S_1)$ and $(B_q)_{q\in Q} \in \mcf_Q(S_2)$, we shall construct 
 $(C_q)_{q\in Q} \in \mcf_Q(S_1 \times S_2)$  such that 
 \[
H^{S_1}_Q(A_q) \cap H^{S_2}_Q(B_q) = H^{S_1 \times S_2}_Q(C_q). 
 \]

Let $C_q = A_q\times B_q \in S_1 \times S_2$ for all $q \in Q$.  If 
\[
A = \bigcap_{q\in Q} A_q \qqand B  = \bigcap_{q\in Q} B_q 
\]
then 
\begin{align*}
C & = \bigcap_{q\in Q} C_q = \bigcap_{q\in Q} (A_q\times  B_q) \\
&  = \bigcap_{q\in Q} A_q \times  \bigcap_{q\in Q} B_q  = A \times  B.
\end{align*}

If $h \in H^{S_1 \times S_2}_Q(C_q)$, then 
\[
C^h = \bigcap_{q\in Q} C_q^h.   
\]
By Lemma~\ref{intersect:lemma:times}, 
\begin{align*} 
A^h \times B^h & = (A \times B)^h = C^h  = \bigcap_{q\in Q} C_q^h \\
& =  \bigcap_{q\in Q} (A_q \times B_q)^h =  \bigcap_{q\in Q} \left( A_q^h \times B_q^h  \right)  \\
& =  \bigcap_{q\in Q}  A_q^h \times  \bigcap_{q\in Q} B_q^h 
\end{align*} 
and so 
\[
A^h = \bigcap_{q\in Q} A_q^h 
\qqand 
B^h = \bigcap_{q\in Q} B_q^h . 
\]
It follows that $h \in H^{S_1}_Q(A_q)$ and $h \in H^{S_2}_Q(B_q)$ 
and so $h \in H^{S_1}_Q(A_q) \cap H^{S_2}_Q(B_q)$ and 
\[
H^{S_1 \times S_2}_Q(C_q) \subseteq H^{S_1}_Q(A_q) \cap H^{S_2}_Q(B_q). 
\]

Conversely, if $h \in  H^{S_1}_Q(A_q) \cap H^{S_2}_Q(B_q)$, 
then 
\[
A^h = \bigcap_{q\in Q} A_q^h \qqand B^h = \bigcap_{q\in Q} B_q^h.  
\]
By Lemma~\ref{intersect:lemma:times}, 
\begin{align*}
C^h & =  (A \times B)^h = A^h \times B^h \\
& = \bigcap_{q\in Q} A_q^h \times \bigcap_{q\in Q} B_q^h 
 = \bigcap_{q\in Q} \left( A_q^h \times \ B_q^h \right) \\ 
& = \bigcap_{q\in Q} \left( A_q  \times \ B_q  \right)^h   = \bigcap_{q\in Q} C_q^h
\end{align*} 
and $h \in H^{S_1 \times S_2}_Q(C_q)$.  
Thus, 
\[
 H^{S_1}_Q(A_q) \cap H^{S_2}_Q(B_q) \subseteq H^{S_1 \times S_2}_Q(C_q). 
\]
It follows that $ H^{S_1}_Q(A_q) \cap H^{S_2}_Q(B_q) = H^{S_1 \times S_2}_Q(C_q)$ 
and so $\mathbf{H}_Q$ 
is closed under intersection.  

To prove that $\mathbf{H}_{\N}^*$ is closed under intersection, it suffices to observe 
that if $(A_q)_{q=1}^{\infty}$ and $(B_q)_{q=1}^{\infty}$ are asymptotically strictly
decreasing sequences of sets, then the sequence  
$(C_q)_{q=1}^{\infty} = (A_q \times B_q)_{q=1}^{\infty} $ 
is also asymptotically strictly decreasing.  
This completes the proof.  
\end{proof}

\bprob 
Are the sets  $\mathbf{H}_Q$ and $\mathbf{H}_{\N}^*$ also closed under unions? 
\eprob

\bprob 
Is there a set $X$ of positive integers with $1 \in X$ such that $X \notin \mathbf{H}_Q$? 
\eprob

\section{$R$-module transformations} 
Let $R$ be a ring, not necessarily commutative, and let $R^{\times}$ be the group of units of $R$. 
Let $M$ be an $R$-module and let $A$ and $A'$ be subsets of $M$.  
We write $A \sim A'$ if there exist 
 $\varepsilon \in R^{\times}$ and $t \in M$ such that 
\[
A' = \varepsilon A + t = \{\varepsilon a + t: a \in A\}.
\]
For example, in  \Z,  we have $A = \{ -7, -5,-4\} \sim A' =  \{0,1,3\}$
because, with $\varepsilon = -1$ and  $t = -4$, there is the relation
$\varepsilon A + t  = (-1)A-4 = \{0,1,3\}= A'$.

Let $M$ be an $R$-module and let $(A_q)_{q \in Q}$ and $(A'_q)_{q \in Q}$ 
be families of subsets of $M$.  We write $(A_q)_{q \in Q} \sim (A'_q)_{q \in Q}$ 
if there exist  $\varepsilon \in R^{\times}$ and $t \in M$ such that  
$A'_q = \varepsilon A_q + t$ for all $q \in Q$.

\bl                                    \label{intersection:lemma:affine} 
The relation $\sim$ is an equivalence relation on the set of subsets of a $R$-module $M$. 
\el

\begin{proof}
Choosing $\varepsilon = 1$ and $t=0$ proves $A \sim A$ for all $A \subseteq M$.

Let $A$ and $A'$ be subsets of $M$ such that $A \sim A'$.   
There exist  $\varepsilon \in R^{\times}$ and $t \in M$ such that  
$A' = \varepsilon A + t$.
Then $\varepsilon^{-1} \in R^{\times}$ and $-\varepsilon^{-1} t  \in M$ and
\begin{align*}
\varepsilon^{-1} A' - \varepsilon^{-1}  t 
& = \varepsilon^{-1} ( \varepsilon A + t ) - \varepsilon ^{-1} t \\
& = \varepsilon^{-1} \varepsilon  A + \varepsilon^{-1}  t - \varepsilon^{-1}  t \\
& =  A
\end{align*}
and so $A' \sim A$. 

Let $A$, $A'$, and $A''$ be subsets of $M$ such that $A \sim A' = \varepsilon A + t$ 
and $A' \sim A'' = \varepsilon' A' + t'$ for some $\varepsilon, \varepsilon' \in R^{\times}$ 
and $t, t' \in M$.  Then 
\begin{align*}
A'' & = \varepsilon' A '+ t' = \varepsilon' (\varepsilon A + t)+ t' \\
& = \varepsilon' \varepsilon A +  \varepsilon'  t + t' \\
& = \varepsilon''  A +  t''
\end{align*}   
with $ \varepsilon''  =  \varepsilon' \varepsilon \in R^{\times}$ and 
$t'' = \varepsilon'  t + t' \in M$, and so $A \sim A''$.
This completes the proof. 
\end{proof}

\bt                     \label{intersection:theorem:H-symmetry} 
Let $M$ be a $R$-module and let $(A_q)_{q\in Q}$ and  $(A'_q)_{q\in Q}$  be families of subsets 
of $M$.  If $(A_q)_{q\in Q} \sim(A'_q)_{q\in Q}$, then $H(A_q) = H(A'_q)$. 
\et

\begin{proof} 
By Lemma~\ref{intersection:lemma:affine}, it suffices to prove that 
$H(A_q) \subseteq \mch(A'_q)$. 

Let $A = \bigcap_{q=1}^{\infty} A_q$ and $A' = \bigcap_{q=1}^{\infty} A'_q$.  
 There exist $\varepsilon \in R^{\times}$ and  $t \in M$ such that 
 \[
 A'_q = \varepsilon A_q + t 
 \]
for all  $q \in Q$. 
If $h \in H(A_q)$, then $hA = \bigcap_{q=1}^{\infty} hA_q$ and  
\begin{align*}
\bigcap_{q=1}^{\infty} hA'_q 
& = \bigcap_{q=1}^{\infty} h\left( \varepsilon A_q + t \right)  = \bigcap_{q=1}^{\infty} \left(  h\varepsilon A_q + ht \right) \\ 
& = \left( \bigcap_{q=1}^{\infty} h \varepsilon A_q \right) + ht  =   \varepsilon \left( \bigcap_{q=1}^{\infty} hA_q \right) + ht \\
& =  \varepsilon hA+ ht  =  h(\varepsilon A + t )\\ 
& = hA'
\end{align*}  
and so $h \in \mch(A'_q)$.  
This completes the proof. 
\end{proof} 

\bprob
Let $(A_q)_{q=1}^{\infty}$ and  $(A'_q)_{q=1}^{\infty}$ be asymptotically strictly decreasing 
sequences of sets of  integers such that 
\[
A'_q = sA_q+t
\]
for all $q \in \N$ and some $s, t \in \Z$ with $s \notin \{ 0, \pm 1\}$.  
Does this imply $H(A_q) = H(A'_q)$?
\eprob

\end{document}